\documentclass[11pt,oneside]{article}
\usepackage[english]{babel}
\usepackage{amssymb}
\usepackage{amsthm}
\usepackage{amsmath}
\usepackage{a4wide}
\usepackage{esvect}
\usepackage{hyperref}
\usepackage{cleveref}
\usepackage{thmtools}
\usepackage{thm-restate}
\usepackage{setspace}
\usepackage{verbatim}
\usepackage[dvips]{graphicx}
\usepackage{t1enc}
\usepackage{color}
\usepackage{pgf,tikz}
\usepackage{mathrsfs}
\usetikzlibrary{arrows}
\usepackage{geometry}
\usepackage{todonotes}
\newtheorem{thm}{Theorem}[section]

\newtheorem{conjecture}[thm]{Conjecture}

\newcommand{\ignore}[1]{}

\title{An improved upper bound on the maximum degree of terminal-pairable complete graphs}
 \author{Ant\'onio Gir\~ao \thanks{Department of Pure Mathematics and Mathematical Statistics, University of Cambridge, Cambridge, UK; \texttt{A.Girao@dpmms.cam.ac.uk}} 
\and
G\'abor M\'esz\'aros \thanks{Department of Mathematical Sciences, The University of Memphis, Memphis, Tennessee;
\texttt{gmszaros@memphis.edu}} }

\begin{document}
\maketitle
\singlespace
\begin{abstract}
A graph $G$ is terminal-pairable with respect to a demand multigraph $D$ on the
same vertex set as $G$, if there exists edge-disjoint paths joining the end vertices of
every demand edge of $D$. In this short note, we improve the upper bound on the largest $\Delta(n)$ with the property that the complete graph on $n$ vertices is terminal-pairable with respect to any demand multigraph of maximum degree at most $\Delta(n)$. This disproves a conjecture originally stated by Csaba, Faudree, Gy\'arf\'as, Lehel and Schelp.
\end{abstract}
\onehalfspace
\section{Introduction} 
The concept of {\itshape terminal-pairability\/} emerged as a practical networking problem and was introduced by Csaba, Faudree, Gy\'arf\'as, Lehel, and Shelp~\cite{CS}. It was further studied by Faudree, Gy\'arf\'as, and Lehel~\cite{mpp,F,pp} and by Kubicka, Kubicki and Lehel~\cite{grid}. Terminal-pairable networks can be defined as follows: given a simple undirected graph $G=(V(G),E(G))$ and an undirected multigraph $D=(V(D), E(D))$  on the same vertex set ($V(D)=V(G))$ we say that $G$ can realize the edges $e_1,\ldots,e_{|E(D)|}$ of $D$ if there exist edge disjoint paths $P_1,\ldots,P_{|E(D)|}$ in $G$ such that $P_i$ joins that endpoints of $e_i$, $i=1,2,\ldots,|E(D)|$. We call $D$ and its edges the {\it demand graph} and the {\it demand edges} of $G$, respectively. Given $G$ and a family $\mathcal{F}$ of (demand)graphs defined on $V(G)$ we call $G$ {\it terminal-pairable} with respect to $\mathcal{F}$ if every demand graph in $\mathcal{F}$ can be realized in $G$.

Given a simple graph $G$, one central question in the topic of terminal-pairability concerns the maximum value of $\Delta$ for which any demand graph $D$ with maximum degree $\Delta(D)\leq \Delta$ can be realized in $G$. Csaba, Faudree, Gy\'arf\'as, Lehel, and Shelp~\cite{CS} studied the above extremal value for complete graphs. Let $K_n^q$ denote the set of demand multigraphs with maximum degree at most $q$ on a complete graph on $n$ vertices. One can easily verify that if $K_n$  is terminal-pairable with respect to $K_n^q$ then $q$ cannot exceed $n/2$. Indeed, consider the demand graph $D$ obtained by replacing every edge in a one-factor by $q$ parallel edges. In order to create edge-disjoint paths routing the endpoints of the demand edges, one needs to use at least two edges in $K_n$, for most of the demand edges. Thus a rather short calculation implies the indicated upper bound on $\Delta(D)$. 

The same authors proved in \cite{CS} that $K_n^q$ can be realized in $K_n$ if $q\leq \frac{n}{4+2\sqrt{3}}$, and conjectured that if ${n\equiv 2\pmod{4}}$, then the upper bound of $n/2$ is attainable, that is, $K_{n}$ is terminal-pairable with respect to $K_n^{n/2}$. This conjecture is also stated in \cite{F}. In this note, we disprove this conjecture by showing that $q/n$ is assymptotically bounded away from $1/2$.

\begin{thm}\label{thm}
If $K_n$ is terminal-pairable with respect to $K_n^{q}$, then $q\leq \frac{13}{27}n + O(1)$.
\end{thm}

We do not know if the newly established upper bound is asymptotically sharp. To this date the best known lower bound on $q$ is $\frac{n}{3} - O(1)$ (see \cite{tpc}).  

\section{Proof of Theorem \ref{thm}}
We may assume $n$ is divisible by $3$. Let $q$ be an even integer. We shall construct a demand graph on $n$ vertices by partitioning the vertex set of $K_n$ into triples, each one forming a triangle where every edge has multiplicity $\frac{q}{2}$. Assume that there exists an edge-disjoint path system $\mathcal{P}$ in $K_n$ that satisfies this demand graph. Note that $e(D)= \frac{nq}{2}$ and at most $n$ demand edges can be realized using exactly one edge of $K_n$ or using $2$ edges within its triple, thus at least $\frac{nq}{2} - n$ demand edges correspond to path of length 2 or more in $\mathcal{P}$. In particular, if $t$ denotes the number of paths of length 2 in $\mathcal{P}$, then the following condition holds due to simple edge counting:
\[2t + 3(\frac{nq}{2} -n-t)\leq \frac{n(n-1)}{2},\]
that is, $t\geq\frac{n}{2}(3q-n-5)$. Hence, if $q$ is sufficiently large (in terms of $n$), then lots of demand edges must be realized through paths of length 2 ("cherries"). 

For a triangle $T_i$, let $\alpha_i$ denote the number of demand edges not in $T_i$ that are resolved in a cherry through any vertex of $T_i$. Also, let $\beta_i$ be the number of demand edges of $T_i$ that are resolved via a cherry with its middle vertex lying outside of $T_i$. Observe that by simple double-counting \[\sum\limits_{i=1}^{\frac{n}{3}}\alpha_i+\beta_i = 2t \geq n(3q-n-5)\]
and therefore there must exist a triangle $T_i$ with $\alpha_i+\beta_i\geq 3(3q-n-5)$. 

Note that between two distinct triangles at most $4$ edges can be solved via paths of length $2$ (every cherry requires two edges between the triangles in $K_n$ and we only have 9 of them). 
This implies that between $T_i$ and any other triangle at most $4 \cdot \left(\frac{n}{3} -1\right)$ demand edges can be solved via cherries. Hence, $4 \cdot \left(\frac{n}{3} -1\right) \geq 3(3q-n-5)$ which implies $q \leq \frac{13}{27}n+1$ as desired.

\section{Additional Remarks}
Terminal-pairability of non-complete graphs has been recently studied by Colucci et al. \cite{btpc}. The authors investigated the extremal value of $q$ for which the set $K_{n,n}^q$ consisting of every demand multigraph on the complete balanced bipartite graph $K_{n,n}$ with maximum degree at most $q$, can be realized. Note that in this variant of the problem the demand graph does not have to be bipartite. It was shown in \cite{btpc} that $K_{n,n}^q$ can be realized in $K_{n,n}$ as long as $q\leq (1-o(1))\frac{n}{4}$. On the other hand, $q$ has to be smaller than $(1+o(1))\frac{n}{3}$ since any one-factor in which every edge has multiplicity $n/3$ is not realizable in $K_{n,n}$. However, we believe that the one-factor construction does not yield a sharp bound:
\begin{conjecture}
There exists $\epsilon >0$ such that if $K_{n,n}^q$ can be realized in $K_{n,n}$, then $q < (\frac{1}{3}-\epsilon)n$, for every $n$ sufficiently large. 
\end{conjecture}

\bibliographystyle{acm}
\bibliography{refs}

\end{document}